\documentclass[preprint,9pt]{elsarticle}
\usepackage{amssymb}
\usepackage{amsmath}
\usepackage{amsthm}
\usepackage{enumerate}

\newtheorem{thrm}{Theorem}[section]

\newtheorem{exam}[thrm]{Example}
\newtheorem{cor}[thrm]{Corollary}

\theoremstyle{definition}
\newtheorem{definition}[thrm]{Definition}
\newtheorem{remark}[thrm]{Remark}

\journal{...}
\input{tcilatex}

\begin{document}

\begin{frontmatter}


\cortext[cor1]{Corresponding author (+903562521616-3087)}

\title{Determinantal and permanental representation of generalized bivariate Fibonacci $p$-polynomials}


\author[rvt]{Kenan KAYGISIZ\corref{cor1}}
\ead{kenan.kaygisiz@gop.edu.tr}
\author[rvt]{Adem SAHIN}
\ead{adem.sahin@gop.edu.tr}
\address[rvt]{Department of Mathematics, Faculty of Arts and Sciences,
Gaziosmanpa\c{s}a University, 60250 Tokat, Turkey}

\begin{abstract}
In this paper, we give some determinantal and permanental
representations of generalized bivariate Fibonacci $p$-polynomials
by using various Hessenberg matrices. The results that we obtained
are important since generalized bivariate Fibonacci $p$-polynomials
are general form of, for example, bivariate Fibonacci and Pell
$p$-polynomials, second kind Chebyshev polynomials, bivariate
Jacobsthal polynomials etc.
\end{abstract}

\end{frontmatter}



\section{Introduction}

Fibonacci numbers and their generalizations have been studying for a long
time. Miles \cite{mil} defined generalized order-$k$ Fibonacci numbers(GO$k$%
F) as,%
\begin{equation}
f_{k,n}=\sum\limits_{j=1}^{k}f_{k,n-j}\   \label{fibmil}
\end{equation}%
for $n>k\geq 2$, with boundary conditions: $f_{k,1}=f_{k,2}=f_{k,3}=\cdots
=f_{k,k-2}=0$ and $f_{k,k-1}=f_{k,k}=1.$\newline

Er \cite{er} defined $k$ sequences of generalized order-$k$ Fibonacci
numbers ($k$SO$k$F) as; for $n>0,$ $1\leq i\leq k$%
\begin{equation*}
f_{k,n}^{\text{ }i}=\sum\limits_{j=1}^{k}c_{j}f_{k,n-j}^{\text{ }i}\ \
\end{equation*}%
with boundary conditions for $1-k\leq n\leq 0,$

\begin{equation*}
f_{k,n}^{\text{ }i}=\left\{
\begin{array}{l}
1\text{ \ \ \ \ \ if \ }i=1-n, \\
0\text{ \ \ \ \ \ otherwise,}%
\end{array}%
\right.
\end{equation*}%
where $c_{j}$ $(1\leq j\leq k)$ are constant coefficients, $f_{k,n}^{\text{ }%
i}$ is the $n$-th term of $i$-th sequence of order-$k$ generalization.

The generalized Fibonacci $p$-numbers \cite{kýlýc} are

\begin{equation}
F_{p}(n)=F_{p}(n-1)+F_{p}(n-p-1)\
\end{equation}%
for $n>p+1$, with boundary conditions $F_{p}(1)=F_{p}(2)=\cdots
=F_{p}(p)=F_{p}(p+1)=1.$\newline

The Fibonacci \cite{alexlupa}, Pell \cite{horadam}, second kind Chebysev
\cite{udrea} and Jacobsthal \cite{tuglu} polynomial are defined as%
\begin{eqnarray*}
f_{n+1}(x) &=&xf_{n}(x)+f_{n-1}(x),\text{ }n\geq 2\text{ with }f_{0}(x)=0,%
\text{ }f_{1}(x)=1 \\
P_{n+1}(x) &=&2xP_{n}(x)+P_{n-1}(x),\text{ }n\geq 2\text{ with }P_{0}(x)=0,%
\text{ }P_{1}(x)=x \\
U_{n+1}(x) &=&xU_{n}(x)-U_{n-1}(x),\text{ }n\geq 2\text{ with }U_{0}(x)=1,%
\text{ }U_{1}(x)=2x \\
J_{n+1}(x) &=&J_{n}(x)+xJ_{n-1}(x),\text{ }n\geq 2\text{ with }J_{0}(x)=0,%
\text{ }J_{1}(x)=1
\end{eqnarray*}%
respectively.

The generalized bivariate Fibonacci $p$-polynomials \cite{tuglu} are%
\begin{equation*}
F_{p,n}(x,y)=F_{p,n-1}(x,y)+F_{p,n-p-1}(x,y)
\end{equation*}%
for $n>p$, with boundary conditions $F_{p,0}(x,y)=0,$ $F_{p,n}(x,y)=x^{n-1},$
$n=1,2,...,p.$\newline

A few terms of $F_{p,n}(x,y)$ for $p=3$ and $p=4$ are

\begin{eqnarray*}
&&0,1,x,x^{2},x^{3},y+x^{4},2xy+x^{5},\ldots \\
&&0,1,x,x^{2},x^{3},x^{4},y+x^{5},2xy+x^{6},\ldots
\end{eqnarray*}%
respectively.

MacHenry \cite{mach1} defined generalized Fibonacci polynomials $%
(F_{k,n}(t)) $ where $t_{i}$ $(1\leq i\leq k)$ are constant coefficients of
the core polynomial%
\begin{equation*}
P(x;t_{1},t_{2},\ldots ,t_{k})=x^{k}-t_{1}x^{k-1}-\cdots -t_{k},
\end{equation*}%
which is denoted by the vector $t=(t_{1},t_{2},\ldots ,t_{k}).$

$F_{k,n}(t)$ is defined inductively by
\begin{eqnarray*}
F_{k,n}(t) &=&0,\text{ }n<1 \\
F_{k,1}(t) &=&1 \\
F_{k,n+1}(t) &=&t_{1}F_{k,n}(t)+\cdots +t_{k}F_{k,n-k+1}(t).
\end{eqnarray*}

\bigskip MacHenry studied on these polynomials and obtained very useful
properties of these polynomials in \cite{mach2, mach3}.

\begin{remark}
\label{remark}Pairs of cognate polynomial sequence \cite{tuglu}%
\begin{equation*}
\begin{tabular}{llll}
\hline
$x$ & $y$ & $p$ & $F_{p,n}(x,y)$ \\ \hline
$x$ & $y$ & $1$ & bivariate Fibonacci polynomials $F_{n}(x,y)$ \\
$x$ & $1$ & $p$ & Fibonacci $p-$polynomials $F_{p,n}(x)$ \\
$x$ & $1$ & $1$ & Fibonacci polynomials $f_{n}(x)$ \\
$1$ & $1$ & $p$ & Fibonacci $p-$numbers $F_{p}(n)$ \\
$1$ & $1$ & $1$ & Fibonacci numbers $F_{n}$ \\
$2x$ & $y$ & $p$ & bivariate Pell $p$-polynomials $F_{p,n}(2x,y)$ \\
$2x$ & $y$ & $1$ & bivariate Pell polynomials $F_{n}(2x,y)$ \\
$2x$ & $1$ & $p$ & Pell $p$-polynomials $P_{p,n}(x)$ \\
$2x$ & $1$ & $1$ & Pell polynomials $P_{n}(x)$ \\
$2$ & $1$ & $1$ & Pell numbers $P_{n}$ \\
$2x$ & $-1$ & $1$ & second kind Chebysev polynomials $U_{n-1}(x)$ \\
$x$ & $2y$ & $p$ & bivariate Jacobsthal $p$-polynomials $F_{p,n}(x,2y)$ \\
$x$ & $2y$ & $1$ & Bivariate Jacobsthal Polynomials $F_{n}(x,2y)$ \\
$1$ & $2y$ & $1$ & Jacobsthal Polynomials $J_{n}(y)$ \\
$1$ & $2$ & $1$ & Jacobsthal Numbers $J_{n}$ \\ \hline
\end{tabular}%
\end{equation*}
\end{remark}

This Remark shows that $F_{p,n}(x,y)$ are general form of all sequences and
polynomials mentioned in Remark \ref{remark}. Therefore, any result obtained
from $F_{p,n}(x,y)$ are valid for all sequences and polynomials mentioned in
Remark \ref{remark}.

Many researchers studied on determinantal and permanental representations of
$k$ sequences of generalized order-$k$ Fibonacci and Lucas numbers. For
example, Minc \cite{min} defined an $n\times n$ (0,1)-matrix $F(n,k),$ and
showed that the permanents of $F(n,k)$ is equal to the generalized order-$k$
Fibonacci numbers (\ref{fibmil}).

The authors (\cite{lee}, \cite{lee2}) defined two $(0,1)$-matrices and
showed that the permanents of these matrices are the generalized Fibonacci (%
\ref{fibmil}) and Lucas numbers. \"{O}cal \cite{oca} gave some determinantal
and permanental representations of $k$-generalized Fibonacci and Lucas
numbers, and obtained Binet's formula for these sequences. K\i l\i c \cite%
{kilic2} gave permanent representation of Fibonacci and Lucas $p$-Numbers.
K\i l\i c \cite{kilic3} studied on permanents and determinants of Hessenberg
matrices. Y\i lmaz and Bozkurt \cite{yil} derived some relationships between
Pell sequences, and permanents and determinants of a type of Hessenberg
matrices.

In this paper, we give some determinantal and permanental representations of
$F_{p,n}(x,y)$ by using various Hessenberg matrices. These results are
general form of determinantal and permanental representations of polynomials
and sequences mentioned Remark \ref{remark}.

\section{The determinantal representations}

In this section, we give some determinantal representations of $F_{p,n}(x,y)
$ by using various Hessenberg matrices. \newline

\bigskip

\begin{definition}
An $n\times n$ matrix $A_{n}=(a_{ij})$ is called lower Hessenberg matrix if $%
a_{ij}=0$ when $j-i>1$ i.e.,%
\begin{equation*}
A_{n}=\left[
\begin{array}{ccccc}
a_{11} & a_{12} & 0 & \cdots & 0 \\
a_{21} & a_{22} & a_{23} & \cdots & 0 \\
a_{31} & a_{32} & a_{33} & \cdots & 0 \\
\vdots & \vdots & \vdots &  & \vdots \\
a_{n-1,1} & a_{n-1,2} & a_{n-1,3} & \cdots & a_{n-1,n} \\
a_{n,1} & a_{n,2} & a_{n,3} & \cdots & a_{n,n}%
\end{array}%
\right]
\end{equation*}
\end{definition}

\begin{thrm}
\label{cahill} \cite{cah} Let $A_{n}$ be an $n\times n$ lower Hessenberg
matrix for all $n\geq 1$ and $\det (A_{0})=1$. Then,
\begin{equation*}
\det (A_{1})=a_{11}
\end{equation*}%
and for $n\geq 2$
\begin{equation*}
\quad \quad \quad \det (A_{n})=a_{n,n}\det (A_{n-1})+\sum\limits_{r=1}^{n-1}
\left[ (-1)^{n-r}a_{n,r}(\prod\limits_{j=r}^{n-1}a_{j,j+1})\det (A_{r-1})%
\right] .
\end{equation*}
\end{thrm}

\bigskip

\begin{thrm}
\label{t1}Let $F_{p,n}(x,y)$ be the generalized bivariate Fibonacci $p$%
-polynomials and $W_{p,n}=(w_{ij})$ an $n\times n$ Hessenberg matrix defined
by
\begin{equation*}
w_{ij}=\left\{
\begin{array}{ll}
i & \text{if \ }i=j-1\text{ }, \\
x & \text{if \ }i=j\text{ }, \\
i^{p}y & \text{if \ }p=i-j,\text{ } \\
0 & \text{otherwise}%
\end{array}%
\right.
\end{equation*}%
that is%
\begin{equation}
W_{p,n}=\left[
\begin{array}{ccccc}
x & i & 0 & \cdots & 0 \\
0 & x & i & \ddots & \vdots \\
\vdots & 0 & x &  & 0 \\
i^{p}y & 0 & \vdots & \cdots &  \\
0 & i^{p}y & 0 &  & 0 \\
\vdots & 0 & \ddots & x & i \\
0 & 0 & \cdots & 0 & x%
\end{array}%
\right] .  \label{kuka}
\end{equation}%
Then%
\begin{equation}
\det (W_{p,n})=F_{p,n+1}(x,y)  \label{tt1}
\end{equation}%
where $n\geq 1$ and $i=\sqrt{-1}.$
\end{thrm}

\bigskip

\begin{proof}
To prove (\ref{tt1}), we use the mathematical induction on $m$. The
result is true for $m=1$ by hypothesis.

Assume that it is true for all positive integers less than or equal
to $m,$ namely $\det (W_{p,m})=F_{p,m}(x,y)$. Then, we have
\begin{eqnarray*}
\det (W_{p,m+1}) &=&q_{m+1,m+1}\det
(W_{p,m})+\sum\limits_{r=1}^{m}\left[
(-1)^{m+1-r}q_{m+1,r}(\prod\limits_{j=r}^{m}q_{j,j+1})\det
(W_{p,r-1})\right]
\\
&=&x\det (W_{p,m})+\sum\limits_{r=1}^{m-p}\left[
(-1)^{m+1-r}q_{m+1,r}(\prod\limits_{j=r}^{m}q_{j,j+1})\det
(W_{p,r-1})\right]
\\
&&+\sum\limits_{r=m-p+1}^{m}\left[
(-1)^{m+1-r}q_{m+1,r}(\prod\limits_{j=r}^{m}q_{j,j+1})\det
(W_{p,r-1})\right]
\\
&=&x\det (W_{p,m})+\left[
(-1)^{p}(i)^{p}y\prod\limits_{j=m-p+1}^{m}i\det (W_{p,m-p})\right]
\\
&=&x\det (W_{p,m})+\left[
(-1)^{p}y(i)^{p}.(i)^{p}\det (W_{p,m-p})\right] \\
&=&x\det (W_{p,m})+y\det (W_{p,m-p})\\
\end{eqnarray*}%
by using Theorem \ref{cahill}. From the hypothesis of induction and
the definition of $F_{p,n}(x,y)$ we obtain
\begin{equation*}
\det (W_{p,m+1})=xF_{p,m+1}(x,y)+yF_{p,m-p+1}(x,y)=F_{p,m+2}(x,y).
\end{equation*}%
Therefore, (\ref{tt1}) holds for all positive integers $n$.
\end{proof}

\bigskip

\begin{exam}
We obtain $6$-th $F_{p,n}(x,y)$ for $p=4$, by using Theorem \ref{t1}%
\begin{equation*}
\ F_{4,6}(x,y)=\det \left[
\begin{array}{ccccc}
x & i & 0 & 0 & 0 \\
0 & x & i & 0 & 0 \\
0 & 0 & x & i & 0 \\
0 & 0 & 0 & x & i \\
i^{4}y & 0 & 0 & 0 & x%
\end{array}%
\right] =y+x^{5}.
\end{equation*}
\end{exam}

\bigskip

\begin{thrm}
\bigskip \label{t2}Let $p\geq 1$ be an integer$,$ $F_{p,n}(x,y)$ be the
generalized bivariate Fibonacci $p$-polynomials and $M_{p,n}=(m_{ij})$ an $%
n\times n$ Hessenberg matrix defined by
\begin{equation*}
m_{ij}=\left\{
\begin{array}{ll}
-1 & \text{if \ }j=i+1, \\
x & \text{if\ \ }i=j, \\
y & \text{if }i-j=p, \\
0\text{\ } & \text{otherwise}%
\end{array}%
\right.
\end{equation*}%
that is
\begin{equation}
M_{p,n}=\left[
\begin{array}{ccccc}
x & -1 & 0 & \cdots & 0 \\
0 & x & -1 & \cdots & 0 \\
0 & 0 & x & \cdots & 0 \\
\vdots & \vdots & \vdots &  & \vdots \\
y & 0 & 0 & \cdots & 0 \\
0 & y & 0 & \cdots & 0 \\
& \vdots & \vdots & \ddots & -1 \\
0 & 0 & \cdots & 0 & x%
\end{array}%
\right] .  \label{beka}
\end{equation}%
Then%
\begin{equation*}
\det (M_{p,n})=F_{p,n+1}(x,y).
\end{equation*}
\end{thrm}

\bigskip

\begin{proof}
Since the proof is similar to the proof of Theorem \ref{t1} by using Theorem %
\ref{cahill} we omit the detail.
\end{proof}

\bigskip

\begin{exam}
\bigskip We obtain $6$-th $F_{p,n}(x,y)$ for $p=3$, by using Theorem \ref{t2}%
\begin{equation*}
F_{3,6}(x,y)=\det \left[
\begin{array}{ccccc}
x & -1 & 0 & 0 & 0 \\
0 & x & -1 & 0 & 0 \\
0 & 0 & x & -1 & 0 \\
y & 0 & 0 & x & -1 \\
0 & y & 0 & 0 & x%
\end{array}%
\right] =2xy+x^{5}.
\end{equation*}
\end{exam}

\section{The permanent representations}

\bigskip

Let $A=(a_{i,j})$ be a square matrix of order $n$ over a ring R. The
permanent of $A$ is defined by%
\begin{equation*}
\text{per}(A)=\sum\limits_{\sigma \in
S_{n}}\prod\limits_{i=1}^{n}a_{i,\sigma (i)}
\end{equation*}%
where $S_{n}$ denotes the symmetric group on $n$ letters.

\begin{thrm}
\label{ocal} \cite{oca} Let $A_{n}$ be an $n\times n$ lower Hessenberg
matrix for all $n\geq 1$ and per$(A_{0})=1.$ Then, $\text{per}(A_{1})=a_{11}$
and for $n\geq 2$
\begin{equation*}
\text{per}(A_{n})=a_{n,n}\text{per}(A_{n-1})+\sum\limits_{r=1}^{n-1}\left[
a_{n,r}(\prod\limits_{j=r}^{n-1}a_{j,j+1})\text{per}(A_{r-1})\right] .
\end{equation*}
\end{thrm}

\bigskip

\begin{thrm}
\bigskip \bigskip \label{t3}Let $p\geq 1$ be an integer$,$ $F_{p,n}(x,y)$ be
the generalized bivariate Fibonacci $p$-polynomials and $H_{p,n}=(h_{rs})$
be an $n\times n$ lower Hessenberg matrix such that%
\begin{equation*}
h_{rs}=\left\{
\begin{array}{ll}
-i & \text{if \ }s-r=1\text{ }, \\
x & \text{if \ }r=s\text{ }, \\
i^{p}y & \text{if \ }p=r-s\text{ } \\
0 & \text{otherwise}%
\end{array}%
\right.
\end{equation*}%
then%
\begin{equation*}
\text{per}(H_{p,n})=F_{p,n+1}(x,y)
\end{equation*}%
where $n\geq 1$ and $i=\sqrt{-1}.$
\end{thrm}

\bigskip

\begin{proof}
This is similar to the proof of Theorem \ref{t1} using Theorem
\ref{ocal}.
\end{proof}

\bigskip

\begin{exam}
\bigskip We obtain $6$-th $F_{p,n}(x,y)$ for $p=3$, by using Theorem \ref{t3}%
\begin{equation*}
F_{3,6}(x,y)=\text{per}\left[
\begin{array}{ccccc}
x & -i & 0 & 0 & 0 \\
0 & x & -i & 0 & 0 \\
0 & 0 & x & -i & 0 \\
-iy & 0 & 0 & x & -i \\
0 & -iy & 0 & 0 & x%
\end{array}%
\right] =\allowbreak 2xy+x^{5}.
\end{equation*}
\end{exam}

\begin{thrm}
\label{t4}Let $p\geq 1$ be an integer$,$ $F_{p,n}(x,y)$ be the
generalized bivariate Fibonacci $p$-polynomials and
$K_{p,n}=(k_{ij})$ be an $n\times n$
lower Hessenberg matrix such that%
\begin{equation*}
k_{ij}=\left\{
\begin{array}{ll}
1 & \text{if \ }j=i+1, \\
x & \text{if \ }i=j, \\
y & \text{if \ }i-j=p, \\
0\text{\ } & \text{otherwise.}%
\end{array}%
\right.
\end{equation*}%
Then%
\begin{equation*}
\text{per}(K_{p,n})=F_{p,n+1}(x,y).
\end{equation*}
\end{thrm}

\bigskip

\begin{proof}
This is similar to the proof of Theorem \ref{t1} by using Theorem
\ref{ocal}.
\end{proof}

\bigskip

We note that the theorems given above are still valid for the sequences and
polynomials mentioned Remark \ref{remark}

\begin{cor}
If we rewrite Theorem \ref{t1}, Theorem \ref{t2}, Theorem \ref{t3} and
Theorem \ref{t4} for $x,y,p$ we obtain the following table.%
\begin{equation*}
\begin{tabular}{llll}
\hline
\textbf{For}$\ \ \mathbf{x}$ & $\mathbf{y}$ & $\mathbf{p}$ & \ \ \ $\det
(W_{p,n})=\det (M_{p,n})=$per$(H_{p,n})=$per$(K_{p,n})=\mathbf{F}_{p,n+1}%
\mathbf{(x,y)},$ \\ \hline
for$\ \ x$ & $y$ & $1$ & \ \ \ $\det (W_{p,n})=\det (M_{p,n})=$per$(H_{p,n})=
$per$(K_{p,n})=\mathbf{F}_{n}\mathbf{(x,y)},$ \\
for $\ x$ & $1$ & $p$ & \ \ \ $\det (W_{p,n})=\det (M_{p,n})=$per$(H_{p,n})=$%
per$(K_{p,n})=$ $\mathbf{F}_{p,n}\mathbf{(x)},$ \\
for$\ \ x$ & $1$ & $1$ & \ \ \ $\det (W_{p,n})=\det (M_{p,n})=$per$(H_{p,n})=
$per$(K_{p,n})=$\textbf{\ }$f_{n}(x),$ \\
for$\ \ 1$ & $1$ & $p$ & \ \ \ $\det (W_{p,n})=\det (M_{p,n})=$per$(H_{p,n})=
$per$(K_{p,n})=$ $\mathbf{F}_{p}\mathbf{(n)},$ \\
for$\ \ 1$ & $1$ & $1$ & \ \ \ $\det (W_{p,n})=\det (M_{p,n})=$per$(H_{p,n})=
$per$(K_{p,n})=$ $\mathbf{F}_{n},$ \\
for$\ \ 2x$ & $y$ & $p$ & \ \ \ $\det (W_{p,n})=\det (M_{p,n})=$\ per$%
(H_{p,n})=$per$(K_{p,n})=$ $\mathbf{F}_{p,n}\mathbf{(2x,y)},$ \\
for$\ \ 2x$ & $y$ & $1$ & \ \ \ $\det (W_{p,n})=\det (M_{p,n})=$per$%
(H_{p,n})=$per$(K_{p,n})=$ $\mathbf{F}_{n}\mathbf{(2x,y)},$ \\
for$\ \ 2x$ & $1$ & $p$ & \ \ \ $\det (W_{p,n})=\det (M_{p,n})=$per$%
(H_{p,n})=$per$(K_{p,n})=$ $\mathbf{P}_{p,n}\mathbf{(x)},$ \\
for$\ \ 2x$ & $1$ & $1$ & \ \ \ $\det (W_{p,n})=\det (M_{p,n})=$per$%
(H_{p,n})=$per$(K_{p,n})=$ $\mathbf{P}_{n}\mathbf{(x)},$ \\
for$\ \ 2$ & $1$ & $1$ & \ \ \ $\det (W_{p,n})=\det (M_{p,n})=$per$(H_{p,n})=
$per$(K_{p,n})=$ $\mathbf{P}_{n},$ \\
for$\ \ 2x$ & $-1$ & $1$ & \ \ \ $\det (W_{p,n})=\det (M_{p,n})=$per$%
(H_{p,n})=$per$(K_{p,n})=\mathbf{U}_{n-1}\mathbf{(x)},$ \\
for$\ \ x$ & $2y$ & $p$ & \ \ \ $\det (W_{p,n})=\det (M_{p,n})=$per$%
(H_{p,n})=$per$(K_{p,n})=$ $\mathbf{F}_{p,n}\mathbf{(x,2y)},$ \\
for$\ \ x$ & $2y$ & $1$ & \ \ \ $\det (W_{p,n})=\det (M_{p,n})=$per$%
(H_{p,n})=$per$(K_{p,n})=$ $\mathbf{F}_{n}\mathbf{(x,2y)},$ \\
for$\ \ 1$ & $2y$ & $1$ & \ \ \ $\det (W_{p,n})=\det (M_{p,n})=$per$%
(H_{p,n})=$per$(K_{p,n})=$ $\mathbf{J}_{n}\mathbf{(y)},$ \\
for$\ \ 1$ & $2$ & $1$ & \ \ \ $\det (W_{p,n})=\det (M_{p,n})=$per$(H_{p,n})=
$per$(K_{p,n})=$\textbf{\ }$J_{n}.$ \\ \hline
\end{tabular}%
\end{equation*}
\end{cor}

\bigskip 

\bigskip

\end{document}